\documentclass{article}

\def\1ox{{ \Omega^1_{\scriptstyle{X}} }}
\def\2ox{{ \Omega^2_{\scriptstyle{X}} }}

\def\ok1{{ \Omega^1_K }}
\def\ok2{{ \Omega^2_K }}

\def\P{{ {\bf P} }}

\def\ra{{ \rightarrow }}

\def\8{{ {\infty } }}

\def\^{{ ^{\wedge} }}

\def\Z{{ {\bf Z } }}

\usepackage{amsmath}
\usepackage{amsfonts}
\usepackage{amscd}
\usepackage{amssymb}

\title{On the topology  of algebraic surfaces and reduction modulo p}
\author{ Dosang Joe and Minhyong Kim}

\begin{document}
\maketitle
\begin{abstract}
We show that the topology of a simply-connected
smooth projective surface is determined by its
algebraic structure modulo p.
\end{abstract}

\section{Observation}
The goal of this note is to make the following simple 
\medskip

{\bf Observation:}
Let $X$ and $Y$ be simply-connected smooth projective surfaces.
Suppose they are isomorphic modulo $p$ for some prime $p \neq 2$
of good reduction.
Then $X$ and $Y$ are homeomorphic by a map preserving the complex orientations.
\medskip

Here, isomorphic modulo $p$ means the following: We can find
a common finitely-generated ring of definition $R$ for
$X$ and $Y$ (that is, such that $X$ and $Y$ can both be
defined by equations with coefficients in $R$), and a
 prime ideal   with residue field $k$
of characteristic
$p$ such that 
$$X\otimes_R \bar{k}\simeq Y\otimes_R \bar{k}$$
That is, $X$ and $Y$ are isomorphic when the coefficients ae
regarded as lying in the algebraic closure of $k$.
`Good reduction' means the two varieties appearing in this isomorphism
are smooth over $\bar{k}$.
For example, they might both be defined over $\Z$, in which case
the meaning of everything is clear.
\medskip

Here is the proof: By Freedman's theorem (\cite{Fr}, theorem 1.5) 
the oriented homeomorphism-type
of a simply-connected oriented smooth 4-manifold $X$ is determined by the
unimodular intersection pairing
$$B: H^2(X,\Z) \otimes H^2(X,\Z) \ra H^4(X,\Z) \simeq \Z$$
on the free abelian group $H^2(X,\Z) $.
We examine this for algebraic surfaces. By the smooth
and proper base-change theorems (\cite{Gr}, chapters 12-14),
$$B_l: [H^2(X,\Z)\otimes \Z_l] \otimes [H^2(X,\Z)\otimes \Z_l] \ra \Z_l,$$
the $\Z_l$-linear extension of $B$, is determined by
$X$ mod $p$ for any $l\neq p$. In particular, the
rank and type (even or odd, by taking $l=2$) of $B$
is determined by reduction mod $p$. On the other hand,
the signature of $X$ is determined by Hirzebruch's signature
theorem (\cite{Hi}, theorem 8.2.2) to be
$$(c_1^2(X)-2c_2(X))/3$$
The numerical invariants in this formula are preserved under
specialization. Therefore,
we conclude that the rank, type, and signature of $B$
are determined by reduction mod $p$.
This concludes the proof when $B$ is indefinite (\cite{Se}, chapter 5,
theorem 6).
But if $B$ is definite, Donaldson's theorem (\cite{Do},\cite{DK} theorem 1.3.1)
says that  $$B \cong \pm (x_1^2+x_2^2+\cdots x_r^2)$$
where $r$ is the rank. But whether or not $B$ is definite (as well
as the sign)
is also determined by  the
rank and signature, and hence, by the reduction modulo $p$. So we are done.

\section{Comments}
Determination of topological invariants of varieties
by modulo $p$ arithmetic
is of course a well-known side-effect of modern arithmetic
geometry.
 
For a smooth and proper variety, 
the Betti numbers, for example, are determined by the
reduction of the variety modulo a (good) prime $p$.

More intricate invariants can also be brought in
(where we always assume that the variety is 
smooth and proper and the prime is
good):

-The Hodge numbers, for example, are determined by reduction modulo
an  {\em ordinary} prime $p$ (\cite{BK} Theorem (0.7)).

-The integral cohomology groups are determined by reduction modulo
{\em two}  primes: this is because two  primes are
sufficient to determine all $H^i(X,\Z_p)$, and then,
$H^i(X,\Z)$ by the universal coefficient theorem.

-For simply-connected varieties, the rational homotopy groups
are determined by reduction modulo  $p$ (\cite{AM}, \cite{DGMS}, \cite{KH}).

-For simply-connected varieties, the integral higher homotopy groups are
 determined by reduction modulo two primes \cite{AM}.

But the case of simply-connected surfaces is
the only one we know of where something as refined
as the homeomorphism-type
is actually determined by reduction modulo $p$.
Our proof is
 of course a consequence of the very powerful classification theorems
for four-manifolds. As such, it appears
 essentially to be
an accident. On the other hand, it would be
interesting to seek out other non-trivial examples of
such theorems, if only to probe their accidental nature.

We point out also that there are parallel results where topological
invariants are preserved by {\em conjugation}, that is,
hitting the coefficients of some defining equations with an automorphism of
the complex numbers.
As is obvious from our proof, this would be true, for example,
for the homeomorphism type of simply-connected smooth projective
surfaces.

One final remark: As is well-known, as a consequence of
the Weil conjectures \cite{De}, the Betti numbers of a smooth proper
variety are determined just by the zeta-function of its reduction mod
$p$. So in this very rigid case of simply-connected surfaces,
it is natural to ask if the homeomorphism type is
determined by the zeta function. However, the simplest possible
case of
$\P^1 \times \P^1$ and $\P^2$ blown up at a point
provides a counter-example. In fact, this  shows that even
the rational homotopy type cannot be determined by the zeta function.
\medskip

{\bf Acknowlegments:}
We are grateful to Igor Dolgachev for interesting conversations
that led to the observation of this paper. We thank
Jong-il Park for suggesting the last counter-example.

M.K.  was supported in part by the National Science Foundation.

{\footnotesize M.K.: DEPARTMENT OF MATHEMATICS, UNIVERSITY OF ARIZONA, TUCSON,
AZ 85721  and KOREA INSTITUTE FOR ADVANCED STUDY, 207-43 CHEONGRYANGRI-DONG
DONGDAEMUN-GU, SEOUL, KOREA, 130-012.  EMAIL:kim@math.arizona.edu

D.J.: DEPARTMENT OF MATHEMATICS, POHANG UNIVERSITY OF SCIENCE AND TECHNOLOGY,
POHANG, KYOUNGBUK, KOREA 790-784. EMAIL: joe@euclid.postech.ac.kr}
 
\end{document}